%&amstex
\input amstex
\input amsppt.sty
\magnification=\magstep1
\hsize=30truecc
\vsize=22.2truecm
\baselineskip=16truept
%\baselineskip=14pt
%\NoBlackBoxes
\nologo
\TagsOnRight
\def\N{\Bbb N}
\def\Z{\Bbb Z}

\def\Q{\Bbb Q}

\def\C{\Bbb C}
\def\l{\left}
\def\r{\right}
\def\bg{\bigg}
\def\({\bg(}
\def\[{\bg[}
\def\){\bg)}
\def\]{\bg]}
\def\t{\text}
\def\f{\frac}

\def\mo{\roman{mod}}

\def\exp{\roman{exp}}
\def\rank{\roman{rank}}
\def\em{\emptyset}
\def\se {\subseteq}

\def\sm{\setminus}

\def\bi{\binom}
\def\eq{\equiv}
\def\cs{\cdots}
\def\ls{\leqslant}
\def\gs{\geqslant}
\def\al{\alpha}

\def\da{\delta}

\def\bi{\binom}
\def\Proof{\noindent{\it Proof}}

\def\Remark{\medskip\noindent{\it  Remark}}

 \topmatter
 \hbox{Israel J. Math. 170(2009), 235--252.}
 \medskip
 \title Zero-sum problems for abelian $p$-groups and covers of the integers by residue classes\endtitle
 \rightheadtext{Zero-sum problems and covers of the integers}
 \author Zhi-Wei Sun\endauthor
 \affil Department of Mathematics, Nanjing University
 \\ Nanjing 210093, People's Republic of China
 \\  zwsun\@nju.edu.cn
 \\ {\tt http://math.nju.edu.cn/$\sim$zwsun}\endaffil
 \medskip
 %\dedicatory In memory of Paul Erd\H os\enddedicatory
\abstract Zero-sum problems for abelian groups and covers of the
integers by residue classes, are two different active topics
initiated by P. Erd\H os more than 40 years ago and investigated by
many researchers separately since then. In an earlier announcement
[Electron. Res. Announc. Amer. Math. Soc. {\bf 9}(2003), 51-60], the
author claimed some surprising connections among these seemingly
unrelated fascinating areas. In this paper we establish further
connections between zero-sum problems for abelian $p$-groups and
covers of the integers. For example, we extend the famous Erd\H
os-Ginzburg-Ziv theorem in the following way: If $\{a_s(\mo\
n_s)\}_{s=1}^k$ covers each integer either exactly $2q-1$ times or
exactly $2q$ times where $q$ is a prime power, then for any
$c_1,\ldots,c_k\in\Z/q\Z$ there exists an $I\se\{1,\ldots,k\}$ such
that $\sum_{s\in I}1/n_s=q$ and $\sum_{s\in I}c_s=0$. Our main
theorem in this paper unifies many results in the two realms and
also implies an extension of the Alon-Friedland-Kalai result on
regular subgraphs.

\endabstract
\thanks 2000 {\it Mathematics Subject Classification}.
Primary 11B75; Secondary 05A05, 05C07, 05E99, 11B25, 11C08, 11D68, 20D60.
\newline\indent Supported by
the National Natural Science Foundation (grant 10871087) of China.
\newline\indent The initial version of this paper was posted
as {\tt arXiv:math.NT/0305369} on May 26, 2003.
\endthanks
\endtopmatter
 \document
\heading {1. Introduction} \endheading

Let $G$ be an abelian group (written additively).
By $\Cal F(G)$ we mean the set of all finite sequences of elements of $G$
(with repetition allowed).
A sequence $\{c_s\}_{s=1}^k\in\Cal F(G)$
(which is often written as $\prod_{s=1}^kc_s$ by A. Geroldinger and his followers)
is called a {\it zero-sum sequence} if $\sum_{s=1}^kc_s=0$.

In 1961 Erd\H os, Ginzburg and Ziv [EGZ] established the following
celebrated theorem which initiated the study of zero-sum sequences.
\proclaim{Theorem 1.1 {\rm (EGZ Theorem)}} Let $n\in\Z^+=\{1,2,3,\ldots\}$.
For any $c_1,\ldots,c_{2n-1}\in\Z$, there
is an $I\se[1,2n-1]=\{1,\ldots,2n-1\}$ with $|I|=n$ such that
$\sum_{s\in I}c_s\eq0\ (\mo\ n)$. In other words, any sequence in $\Cal F(\Z_n)$
of length $2n-1$ contains a zero-sum subsequence of length $n$, where $\Z_n=\Z/n\Z$ is the additive group
of residue classes modulo $n$.
\endproclaim

The EGZ theorem can be easily reduced to the case where $n$ is a
prime (and hence $\Z_n$ is a field), and then deduced from the
well-known Cauchy-Davenport theorem or the Chevalley-Warning
theorem. (See, e.g., Nathanson [N, pp.\,48--51], and Geroldinger and Halter-Koch [GH, p.\,349].) It remains valid
if we replace the cyclic group $\Z_n$ by an arbitrary abelian
group of order $n$. (Cf. T. Tao and V. Vu [TV, pp.\,350--351].)

Let $G$ be a finite abelian group. When $|G|>1$, there is a unique sequence
$d_1,\ldots,d_r$ of positive integers with $d_1>1$ and $d_i\mid d_{i+1}$ for $1\ls i<r$
such that $G$ is isomorphic to the direct sum
$$\Z_{d_1}\oplus\Z_{d_2}\oplus\cdots\oplus\Z_{d_r};$$
in this case, $r$ is $\rank(G)$ (the rank of $G$) and $d_r$ is
$\exp(G)$ (the exponent of $G$), and we define
$$d^*(G)=\sum_{i=1}^r(d_i-1).\tag1.1$$
If $|G|=1$, then $\rank(G)=\exp(G)=1$ and we set $d^*(G)=0$.
Clearly $d^*(G)+1\ls|G|$.

Let $G$ be a finite abelian group written additively. By $s(G)$ we denote the smallest
positive integer $k$ such that any sequence in $\Cal F(G)$ of length $k$ has a zero-sum subsequence of length $\exp(G)$.
For any $n\in\Z^+$, we have $s(\Z_n)=2n-1$ by the EGZ theorem, and $s(\Z_n\oplus\Z_n)=4n-3$
by the Kemnitz conjecture  proved by C. Reiher [Re] (see also [SC]),
and $s(\Z_d\oplus\Z_n)=2(d+n)-3$ by Theorem 5.8.3 of Geroldinger and Halter-Koch [GH, p.\,362]
where $d$ is any positive divisor of $n$.
The reader is referred to the survey [GG06], and the recent papers [E] and [EEGKR]
for various problems and results on $s(G)$.

Let $G$ be an abelian group of order $n$. For any $\{c_s\}_{s=1}^n\in\Cal F(G)$,
as the following elements
$$0,\ c_1,\ c_1+c_2,\ \ldots,\ c_1+c_2+\cdots+c_{n}$$ cannot be distinct by the pigeon-hole principle,
we have $\sum_{s\in I}c_s=0$ for some $\em\not=I\se[1,n]$; furthermore, $\{c_s\}_{s=1}^n$
has a zero-sum subsequence $\{c_s\}_{s\in I}$ with $\em\not=I\se[1,n]$ and
$\sum_{s\in I}1/\t{ord}(c_s)\ls1$,
by a celebrated theorem of Geroldinger [G93] (which was re-proved later by Elledge and Hurlbert [EH]
via graph pebbling). The {\it Davenport constant} $D(G)$ of  $G$
is defined as the smallest positive
integer $k$ such that any sequence $\{c_s\}_{s=1}^k\in\Cal F(G)$ has a nonemepty zero-sum subsequence.
(Note that we essentially impose no restriction on the length of the required zero-sum subsequence.)
By the above, $D(G)\ls n=|G|$.
In 1966
Davenport showed that if $K$ is an algebraic number field with
ideal class group $G$, then $D(G)$ is the maximal number of prime
ideals (counting multiplicity) in the decomposition of an
irreducible algebraic integer in $K$. The reader may consult
Theorem 5.1.5 of [GH, pp.\,305--306] for further results in this direction.

It is easy to see that $D(\Z_n)=n$ for any $n\in\Z^+$.
For an abelian $p$-group $G$ with $p$ a prime,
 $D(G)$ is greater than $d^*(G)$ by a constructive example;
on the other hand, in 1969 Olson [O] used the knowledge of group rings to show
that $D(G)\ls d^*(G)+1$ (and hence $D(G)=d^*(G)$+1).
Olson's original method has been further refined and explored
by many researchers, see, e.g., [GGH].

\proclaim{Theorem 1.2 {\rm (Olson's Theorem)}} Let $p$ be a prime and let $G$
be an additive abelian $p$-group.
Given $c,c_1,\ldots,c_{d^*(G)+1}\in G$ we have
$$\sum\Sb I\se[1,d^*(G)+1]\\ \sum_{s\in I}c_s=c\endSb(-1)^{|I|}\eq0\ \ (\mo\ p),$$
and in particular $\{c_s\}_{s=1}^{d^*(G)+1}$ has a nonempty zero-sum subsequence.
\endproclaim

Let $q$ be a power of a prime $p$, and let $c\in\Z_q$ and $\{c_s\}_{s=1}^{2q-1}\in\Cal F(\Z_q)$.
By Olson's theorem in the case
$G=\Z_q^2=\Z_q\oplus\Z_q$, we have
$$\sum\Sb I\se[1,d^*(\Z_q^2)+1]\\q\mid|I|,
 \ \sum_{s\in I}c_s=c\endSb(-1)^{|I|}\eq0\ \ (\mo\ p).$$
 In other words,
 $$\bg|\bg\{I\se[1,2q-1]:\, |I|=q\ \t{and}\ \sum_{s\in I}c_s=c\bg\}\bg|
  \eq [\![c=0]\!]\ \ (\mo\ p),$$
 where for a predicate $P$ we let $[\![P]\!]$ be $1$ or $0$
 according as $P$ holds or not.
Thus, Olson's theorem implies the EGZ theorem.

Let $G$ be a finite abelian group. A zero-sum sequence $\{c_s\}_{s=1}^k$
is called a {\it minimal zero-sum sequence} if $\sum_{s\in I}c_s=0$
for no $\em\not=I\subset [1,k]$. Though we don't study minimal zero-sum sequences in this paper,
the reader is still recommended to see [GG99], [GGS] and [LS] for some results on
minimal zero-sum sequences.

 Now we turn to covers of the integers by residue classes.

 For $a\in\Z$ and $n\in\Z^+$ we call
 $$a(n)=a+n\Z=\{a+nx:\, x\in\Z\}$$
 a residue class with modulus $n$. For a finite system
 $$A=\{a_s(n_s)\}_{s=1}^k\tag1.2$$
 of residue classes, its {\it covering function}
 $$w_A(x)=|\{1\ls s\ls k:\, x\in a_s(n_s)\}|$$
 is periodic modulo the least
 common multiple $N_A=[n_1,\ldots,n_k]$ of the moduli $n_1,\ldots,n_k$. Sun [S97, S99]
 called $m(A)=\min_{x\in\Z}w_A(x)$
 the {\it covering multiplicity} of (1.2).
 One can easily verify the following basic property:
 $$\sum_{s=1}^k\f1{n_s}=\f1{N_A}\sum_{x=0}^{N_A-1}w_A(x)\gs m(A).\tag1.3$$
 Further properties of the covering function $w_A(x)$ can be found in [S03a, S04].

  If $\bigcup_{s=1}^ka_s(n_s)=\Z$ (i.e., $m(A)\gs1$), then
 we call (1.2) a {\it cover} (or {\it covering system}) of $\Z$.
 This concept was first introduced by Erd\H os in
 the early 1930's (cf. [E50]), and many surprising applications have been found
 (cf. [Cr], [Gra], [S00] and [S01]).
 Erd\H os was very proud of this invention; in [E97] he said:
 {\it ``Perhaps my favorite problem of all concerns covering systems}."

 For $m\in\Z^+$, if $m(A)\gs m$ then $A$ is said to be an {\it $m$-cover}
 of $\Z$; general $m$-covers were first studied by the author in [S95].
 It is easy to construct an $m$-cover of $\Z$
 which cannot be split into two covers of $\Z$ (cf. [PS, Example 1.1]).

 If $w_A(x)=m$ for all $x\in\Z$, then we call $(1.2)$ an {\it exact $m$-cover} of $\Z$.
 (Note that in this case we have $\sum_{s=1}^k1/n_s=m$ by (1.3).)
 Clearly $m$ copies of $0(1)$ form a trivial exact $m$-cover of $\Z$.
 Using a graph-theoretic argument
 Zhang [Z91] proved that for each $m=2,3,\ldots$ there are infinitely many exact $m$-covers of
 $\Z$ which cannot split into an exact $n$-cover and an exact $(m-n)$-cover with
 $0<n<m$; such an exact $m$-cover is said to be {\it irreducible}. In 1973
 Choi supplied the following example of an irreducible exact
 $2$-cover of $\Z$:
 $$\{1(2);0(3);2(6);0,4,6,8(10);
 1,2,4,7,10,13(15);5,11,12,22,23,29(30)\}.$$
 Zhang [Z91] showed that the residue classes
 $$1(2),0(6),0(10),0(14),1(15),1(15),8(21),8(21),36(105),36(105),$$
 together with some residue classes modulo 210,
 form an irreducible exact $3$-cover of $\Z$.
 In 1992 Sun [S92] proved that if $\{a_s(n_s)\}_{s=1}^k$ forms an exact $m$-cover
 of $\Z$ then for each $n=0,1,\ldots,m$ there are at least $\bi mn$
 subsets $I$ of $[1,k]$ with $\sum_{s\in I}1/n_s=n$.

 There are many problems and results on covers of $\Z$; the reader
 may consult sections F13 and F14 of the book [Gu, pp.\,383--390],
 the survey [P-S], and the recent papers [S05a] and [FFKPY].

 Now we mention some properties of covers of $\Z$ related to Egyptian
 fractions. The first nontrivial result of this nature is
 the following one discovered by Zhang [Z89]
 with the help of the Riemann zeta function: If $\{a_s(n_s)\}_{s=1}^k$
 forms a cover of $\Z$, then
 $$\sum_{s\in I}\f1{n_s}\in\Z
 \quad \t{for some nonempty}\ I\se[1,k].\tag1.4$$
 The following theorem contains two different extensions of Zhang's result.

 \proclaim{Theorem 1.3} Let $\{a_s(n_s)\}_{s=1}^k$ be an $m$-cover of $\Z$, and let $m_1,\ldots,m_k\in\Z^+$.

 {\rm (i) (Sun [S95, S96])} There are at least $m$ positive
integers in the form
 $\sum_{s\in I}m_s/n_s$ with $I\se[1,k]$.

 {\rm (ii) (Pan and Sun [PS])} For any $J\se[1,k]$ there are at
 least $2^m$ subsets $I$ of $[1,k]$ with $\sum_{s\in
 I}m_s/n_s-\sum_{s\in J}m_s/n_s\in\Z$.
 \endproclaim

 Note that a residue class $a(n)=a+n\Z$ is a coset of the subgroup $n\Z$ of the additive group $\Z$.
 There are also some investigations on covers of a general group by left cosets of subgroups,
 see, e.g., [S06] and the references therein. Gao and Geroldinger [GG03] reduced
 some zero-sum problems to the study of covering a certain subset of an abelian group
 by few proper cosets. However, in this paper we are only interested in covers of the integers
 and their surprising connections with zero-sum problems.

  The purpose of this paper is to show that some classical results of zero-sum nature,
 such as the EGZ theorem and Olson's theorem, are special cases
 of our general results on covers of $\Z$.
 The key point is to compare Davenport constants of abelian $p$-groups with
covering multiplicities of covers of $\Z$.

   In Section 2 we state our main results connecting zero-sum problems for abelian $p$-groups
with covers of the integers. A more general theorem will be presented in the third section,
together with some consequences; its proof will be given in Section 4.

 \heading{2. Connections between Zero-sum Sequences and Covers of $\Z$}\endheading

 The following theorem reveals unexpected connections between zero-sum sequences and covers of $\Z$.

\proclaim{Theorem 2.1 {\rm (Main Theorem)}} Let
 $G$ be an additive abelian $p$-group where $p$ is a prime.
 Suppose that $A=\{a_s(n_s)\}_{s=1}^k$ is a $(d^*(G)+p^h)$-cover of $\Z$
 with $h\in\N=\{0,1,\ldots\}$. Let $c_1,\ldots,c_k\in G$ and $m_1,\ldots,m_k\in\Z$. Then
$$\bg|\bg\{I\se[1,k]:\, \sum_{s\in I}c_s=c\ \t{and}\
\sum_{s\in I}\f{m_s}{n_s}\in\al+p^h\Z\bg\}\bg|\not=1\tag2.1$$
for any $c\in G$ and rational number $\al$.
In particular, $\{c_s\}_{s=1}^k$ has a zero-sum subsequence $\{c_s\}_{s\in I}$
with $\em\not=I\se[1,k]$ satisfying the restriction $\sum_{s\in I}m_s/n_s\in p^h\Z$.
 \endproclaim

 Now we deduce various
 consequences from Theorem 2.1.

 \proclaim{Corollary 2.1} Let $A=\{a_s(n_s)\}_{s=1}^k$
 be a $(d^*(G)+1)$-cover of $\Z$ where $G$ is an additive abelian $p$-group.
 Let $m_1,\ldots,m_k\in\Z$. Then any sequence $\{c_i\}_{i=1}^k\in \Cal F(G)$ has a zero-sum subsequence
 $\{c_s\}_{s\in I}$ with $\em\not=I\se[1,k]$ such that $\sum_{s\in I}m_s/n_s\in\Z$.
 \endproclaim

 \Proof. Simply apply Theorem 2.1 with $h=0$. \qed

 \Remark\ 2.1. (a) For an abelian $p$-group $G$,
  if we apply Corollary 2.1 to the trivial $(d^*(G)+1)$-cover
  consisting of $d^*(G)+1$ copies of $0(1)$, we then
  obtain Olson's result $D(G)\ls d^*(G)+1$. We conjecture that
  the desired result in Corollary 2.1 still holds
  when $G$ is a finite abelian group and $A=\{a_s(n_s)\}_{s=1}^k$ is a $D(G)$-cover of $\Z$.

  (b) That (1.4) holds for any cover $(1.2)$, follows from Corollary 2.1 in the case $G=\{0\}$.

  (c) If we apply Corollary 2.1 to the trivial group $G=\{0\}$ and the trivial cover $\{r(n)\}_{r=0}^{n-1}$
  (where $n\in\Z^+$), then we get the basic result $D(\Z_n)=n$ in zero-sum theory.

 \proclaim{Corollary 2.2} Let $\{a_s(n_s)\}_{s=1}^k$ be an $m$-cover of $\Z$ and let $m_1,\ldots,m_k\in\Z$.
 Assume that $m$ is a prime power $($i.e., $m=p^n$ for some prime $p$ and $n\in\N)$.
 Then $\sum_{s\in I}m_s/n_s\in m\Z$ for some $\em\not=I\se[1,k]$, and in particular
 $$\sum_{s\in I}\f1{n_s}\in m\Z^+\quad\t{for some}\ I\se[1,k].\tag2.2$$
 Moreover, for any $J\se[1,k]$
 there is an $I\se[1,k]$ with $I\not=J$ such that
 $$\sum_{s\in I}\f{m_s}{n_s}-\sum_{s\in J}\f{m_s}{n_s}\in m\Z.$$
 \endproclaim
 \Proof. Just apply Theorem 2.1 with $G=\{0\}$ and $\al=\sum_{s\in J}m_s/n_s$. \qed

 \Remark\ 2.2. Corollary 2.2 is a new extension
 of Zhang's result that (1.4) holds if (1.2) is a cover of $\Z$.
 We conjecture that the corollary remains valid if we remove the condition that
 $m$ is a prime power. In the special case $n_1=\cdots=n_k=1$, this conjecture
 yields the basic fact $D(\Z_m)=m$.
 \medskip

 Theorem 2.1 also implies the following central result of this paper.

 \proclaim{Theorem 2.2} Let $G$ be an abelian $p$-group with $p$ a prime, and let $q>d^*(G)$ be a power of $p$
 $($e.g., $q=|G|)$.

 {\rm (i)} Let $A=\{a_s(n_s)\}_{s=1}^k$
 with $\{w_A(x):\,x\in\Z\}\se[d^*(G)+q,2q]$.
 Then any $\{c_s\}_{s=1}^k\in\Cal F(G)$ has a zero-sum subsequence
 $\{c_s\}_{s\in I}$ with $I\se[1,k]$
 and $\sum_{s\in I}1/n_s=q$.

 {\rm (ii)} Let $A=\{a_s(n_s)\}_{s=1}^k$ be an exact $3q$-cover of $\Z$.
 Then any zero-sum sequence  $\{c_s\}_{s=1}^k\in\Cal F(G\oplus G)$
 has a zero-sum subsequence  $\{c_s\}_{s\in I}$ with $I\se[1,k]$
 and $\sum_{s\in I}1/n_s=q$.
\endproclaim
\Proof. (i) By Theorem 2.1, $\{c_s\}_{s=1}^k$ has a zero-sum subsequence
$\{c_s\}_{s\in I}$ with $\em\not=I\se[1,k]$ satisfying $\sum_{s\in I}1/n_s\in q\Z$.
Observe that
$$\sum_{s\in I}\f1{n_s}\ls\sum_{s=1}^k\f1{n_s}=\f1{N_A}\sum_{x=0}^{N_A-1}w_A(x)\ls 2q.$$
If $\sum_{s\in I}1/n_s\not=q$, then $\sum_{s\in I}1/n_s=2q$ and $w_A(x)=2q$ for all $x\in\Z$.
When $A$ is an exact $2q$-cover of $\Z$, for the system $A_*=\{a_s(n_s)\}_{s=1}^{k-1}$
we still have $\{w_{A_*}(x):\,x\in\Z\}\se[d^*(G)+q,2q]$, hence
$\{c_s\}_{s=1}^k$ has a zero-sum subsequence $\{c_s\}_{s\in I_*}$ with
$\em\not=I_*\se[1,k-1]$ and $\sum_{s\in I_*}1/n_s=q$.
This completes the proof of part (i).

(ii) Note that $d^*(G\oplus G)=2d^*(G)$.
As $3q-1>d^*(G\oplus G)+q$, $A_*=\{a_s(n_s)\}_{s=1}^{k-1}$ is a $(d^*(G\oplus G)+q)$-cover of $\Z$.
 Applying Theorem 2.1 to the system $A_*$
 we find that $\{c_s\}_{s=1}^k$ has a zero-sum subsequence
 $\{c_s\}_{s\in I}$ with  $\em\not=I\se[1,k-1]$ and $n=\sum_{s\in I}1/n_s\in q\Z$.
 As $n<\sum_{s=1}^k1/n_s=3q$, $n$ is $q$ or $2q$. If $n=2q$,
 then for $\bar I=[1,k]\sm I$ we have
 $$\sum_{s\in \bar I}c_s=\sum_{s=1}^kc_s-\sum_{s\in I}c_s=0
 \ \t{and}\ \sum_{s\in\bar I}\f1{n_s}
 =\sum_{s=1}^k\f1{n_s}-n=3q-2q=q.$$
 This concludes the proof. \qed

\Remark\ 2.3. It is interesting to view $1/n_s$ in Theorem 2.2 as a weight of $s\in[1,k]$.
In the case $n_1=\cdots=n_k=1$, part (i) yields the EGZ theorem for abelian $p$-groups,
and part (ii) gives
 Lemma 3.2 of Alon and Dubiner [AD], which is an indispensable lemma
 in the study of the Kemnitz conjecture (cf. [Ro] and [Re]).
 Note that our Theorem 2.2(i) is quite different from the so-called weighted EGZ theorem
 proved by Grynkiewicz [Gry].
\medskip

Theorem 2.2 tells that our following conjecture holds when $n$ is a prime power.

 \proclaim{Conjecture 2.1} Let $G$ be a finite abelian group of order $n$.

 {\rm (i)} If $\{a_s(n_s)\}_{s=1}^k$ covers each integer either exactly $2n-1$ times or exactly $2n$ times,
 then any $\{c_s\}_{s=1}^k\in\Cal F(G)$ has a zero-sum subsequence $\{c_s\}_{s\in I}$ with
 $I\se[1,k]$ and $\sum_{s\in I}1/n_s=n$.

 {\rm (ii)} When $\{a_s(n_s)\}_{s=1}^k$ forms an exact $3n$-cover of $\Z$,
 any zero-sum sequence $\{c_s\}_{s=1}^k\in\Cal F(G\oplus G)$
has a zero-sum subsequence $\{c_s\}_{s\in I}$ with $I\se[1,k]$
and $\sum_{s\in I}1/n_s=n$.
 \endproclaim

 An undirected graph is said to be {\it $q$-regular} if all the vertices have degree $q$.
 In 1984 Alon, Friedland and Kalai [AFK1, AFK2]
 proved that if $q$ is a prime power then
 any loopless (undirected) graph with average degree bigger than $2q-2$
 and maximum degree at most $2q-1$ must contain a $q$-regular
 subgraph. Now we apply Theorem 2.1 to
 strengthen this result.

\proclaim{Theorem 2.3} Let $G$ be a loopless graph of $l$ vertices
with the edge set $\{1,\ldots,k\}$. Suppose that all the vertices
of $G$ have degree not exceeding $2p^{n}-1$ and that $\{a_s(n_s)\}_{s=1}^k$
forms an $(l(p^{n}-1)+p^h)$-cover of $\Z$, where $p$ is a prime
and $n,h\in\N$.
Then, for any $m_1,\ldots,m_k\in\Z$, there exists
a $p^n$-regular subgraph $H$ of $G$ with
$\sum_{s\in E(H)}m_s/n_s\in p^h\Z$,
where $E(H)$ denotes the edge set of $H$.
\endproclaim

\Proof. Let $v_1,\ldots,v_l$ be all the vertices of graph $G$.
For $s\in[1,k]$ and $t\in[1,l]$, set
$$\da_{st}=[\![v_t\ \t{is an endvertex of the edge}\ s]\!].$$
Note that $\sum_{s=1}^k\da_{st}$
is just the degree $d_G(v_t)$.

By Theorem 2.1 in the case $G=\Z_{p^n}^l$,
there is a nonempty $I\se[1,k]$
such that $\sum_{s\in I}m_s/n_s\in p^h\Z$ and
$p^{n}\mid\sum_{s\in I}\da_{st}$ for all $t\in[1,l]$.
 Let $V_I$ be the set of
vertices incident with edges in $I$, and let $H$ be the subgraph
$(V_I,I)$ of graph $G$. As $\sum_{s\in I}\da_{st}\ls
d_G(v_t)\ls 2p^{n}-1$, $p^{n}\mid\sum_{s\in I}\da_{st}$ for all $t\in[1,l]$
if and only if $d_H(v)=p^{n}$ for all $v\in V_I$ (i.e., $H$ is a
$p^n$-regular subgraph of $G$). This concludes the proof. \qed

\Remark\ 2.4. For the graph $G$ in Theorem 2.3, clearly
$$k\gs l(p^{n}-1)+1\Leftrightarrow 2k>l(2p^{n}-2)\Leftrightarrow
\sum_{v\in V(G)}d_G(v)>(2p^{n}-2)|V(G)|$$
where $V(G)$ is the vertex set of graph $G$.
So Theorem 2.3 in the case $h=0$ and $n_1=\cs=n_k=1$ implies
the Alon-Friedland-Kalai result.

\heading{3. A General Theorem and its Consequences}\endheading

 Let $\Omega$ be the ring of all algebraic integers. For
 $\omega_1,\omega_2,\gamma\in\Omega$, by $\omega_1\eq\omega_2\ (\mo\ \gamma)$
 we mean $\omega_1-\omega_2\in\gamma\Omega$.
 For $a,b,m\in\Z$ it is well known that $a-b\in m\Omega$ if and
 only if $a-b\in m\Z$ (see, e.g. [IR, p.\,68]).
 For $m\in\Z$ and a root $\zeta$ of unity,
 if $\zeta\eq0\ (\mo\ m)$ then $1=\zeta\zeta^{-1}\eq0\ (\mo\ m)$
 (since $\zeta^{-1}\in\Omega$)
 and hence $m$ must be $1$ or $-1$.

  Theorem 1.2 of zero-sum nature, Theorem 1.3(i) on covers of $\Z$,
  and our useful Theorem 2.1 are special cases
  of the following general theorem (which is inevitably complicated since it
  unifies many results).

  \proclaim{Theorem 3.1} Let
 $G$ be an additive abelian $p$-group where $p$ is a prime.
 Suppose that $A=\{a_s(n_s)\}_{s=1}^k$ is a $(d^*(G)+p^h+\Delta)$-cover of $\Z$
 with $h,\Delta\in\N$. Let $m_1,\ldots,m_k\in\Z$
 and $c,c_1,\ldots,c_k\in G$. Let $\al$ belong to the rational field $\Q$, and set
$$\Cal I=\bg\{I\se[1,k]:\, \sum_{s\in I}c_s=c\ \t{and}\
\sum_{s\in I}\f{m_s}{n_s}\in\al+p^h\Z\bg\}.\tag3.1$$ Let
$P(x_1,\ldots,x_k)\in\Q[x_1,\ldots,x_k]$ have degree not exceeding
$d\in\Z^+$ and $P([\![1\in I]\!],\ldots,[\![k\in I]\!])\in\Z$ for all
$I\se[1,k]$ with $\sum_{s\in I}m_s/n_s-\al\in\Z$. Then, either we have the inequality
$$|\{P([\![1\in I]\!],\ldots,[\![k\in I]\!])\ \mo\ p:\, I\in\Cal I\}|>1+\f{\Delta}d,\tag3.2$$
or $|\{I\in\Cal I:\, P([\![1\in I]\!],\ldots,[\![k\in I]\!])\in p\Z\}|\not=1$ and furthermore
$$\sum\Sb I\in\Cal I\\p\mid P([\![1\in I]\!],\ldots,[\![k\in I]\!])\endSb
(-1)^{|I|} e^{2\pi i\sum_{s\in I}a_sm_s/n_s}\eq0\ \ (\mo\ p).\tag3.3$$
 \endproclaim

\Remark\ 3.1. (a) By taking $P(x_1,\ldots,x_k)=0$ in Theorem 3.1, we get the congruence
$$\sum_{I\in\Cal I}(-1)^{|I|}e^{2\pi i\sum_{s\in I}a_sm_s/n_s}\eq0\ (\mo\ p)$$
under the conditions of Theorem 3.1. (Thus Theorem 2.1 follows from Theorem 3.1.)
In the case $h=0$ and $n_1=\cdots=n_k=1$, this yields Theorem 1.2 of Olson.

(b) When $G$ is an elementary abelian $p$-group, and $h=0$ and $d=1$,
Theorem 3.1 is equivalent to the first part of the Main Theorem
in the announcement [S03b].

 \proclaim{Corollary 3.1} Let $\{a_s(n_s)\}_{s=1}^k$
 be an $m$-cover of $\Z$. Let $m_1,\ldots,m_k\in\Z$, and let
 $\mu_1,\ldots,\mu_k$ be rational numbers such that
 $\sum_{s\in I}\mu_s\in\Z$ for all those $I\se[1,k]$ with $\sum_{s\in I}m_s/n_s\in\Z$.
 For any prime $p$, if there is no $\em\not=I\se[1,k]$ such that $\sum_{s\in I}m_s/n_s\in\Z$ and
 $\sum_{s\in I}\mu_s\in p\Z$, then
 $$\bg|\bg\{\sum_{s\in I}\mu_s\ \mo\ p:\, \em\not=I\se[1,k]\ \t{and}
 \ \sum_{s\in I}\f{m_s}{n_s}\in\Z\bg\}\bg|\gs m.\tag3.4$$
 \endproclaim
 \Proof. Apply Theorem 3.1 with $G=\{0\}$, $h=\al=0$ and $P(x_1,\ldots,x_k)=\sum_{s=1}^k\mu_sx_s$. \qed

 \Remark\ 3.2. (a) Under the conditions of Corollary 3.1, if $\mu_1,\ldots,\mu_k$ are positive then
 by taking a prime $p>\mu_1+\cdots+\mu_k$ we get that
 $$\bg|\bg\{\sum_{s\in I}\mu_s:\, \em\not=I\se[1,k]\ \t{and}
 \ \sum_{s\in I}\f{m_s}{n_s}\in\Z\bg\}\bg|\gs m;$$
 in particular,
$$\bg|\bg\{|I|:\, \em\not=I\se[1,k]\ \t{and}
 \ \sum_{s\in I}\f{m_s}{n_s}\in\Z\bg\}\bg|\gs m.\tag3.5$$
Theorem 1.3(i) follows if we set $\mu_s=m_s/n_s$ for $s=1,\ldots,k$.

  (b) In the special case $n_1=\cdots=n_k=1$, Corollary 3.1 gives the following result:
If $c_1,\ldots,c_k\in\Z_p$ with $p$ a prime, and $\sum_{s\in I}c_s=0$ for no $\em\not=I\se[1,k]$,
then $|\{\sum_{s\in I}c_s:\, \em\not= I\se[1,k]\}|\gs k$.
\medskip

 \proclaim{Corollary 3.2} Let $\{a_s(n_s)\}_{s=1}^k$ be an $m$-cover of $\Z$.
 Given $m_1,\ldots,m_k\in\Z$ and $J\se[1,k]$,
 we have
 $$\bg|\bg\{\sum_{s\in I}\f{a_sm_s}{n_s}-\f{|I|}2:\,
 I\se[1,k]\ \t{and}\ \sum_{s\in I}\f{m_s}{n_s}-\sum_{s\in J}\f{m_s}{n_s}\in\Z\bg\}\bg|>m.\tag3.6$$
 \endproclaim
 \Proof. Let $N$ be the least common multiple of $2,n_1,\ldots,n_k$. Set
 $$\Cal I=\bg\{I\se[1,k]:\, \sum_{s\in I}\f{m_s}{n_s}-\sum_{s\in J}\f{m_s}{n_s}\in\Z\bg\}
 \ \t{and}\ l=\max_{I\se[1,k]}\bg|N\sum_{s\in I}\mu_s\bg|$$
 where $\mu_s=a_sm_s/n_s-1/2$. Choose a prime $p>\max\{|\Cal I|,2l\}$ and set
 $$P(x_1,\ldots,x_k)=N\(\sum_{s=1}^k\mu_sx_s-\sum_{s\in J}\mu_s\)\in\Z[x_1,\ldots,x_k].$$
Since $l<p/2$, for $I_1,I_2\se[1,k]$ we have
$$N\sum_{s\in I_1}\mu_s\eq N\sum_{s\in I_2}\mu_s\ (\mo\ p)\iff \sum_{s\in I_1}\mu_s=\sum_{s\in I_2}\mu_s.$$

  Now assume that (3.6) fails. Then
  $$|\{P([\![1\in I]\!],\ldots,[\![k\in I]\!])\ \mo\ p:\, I\in\Cal I\}|\ls m.$$
  Applying Theorem 3.1 with $G=\{0\}$, $h=0$ and $\al=\sum_{s\in J}m_s/n_s$,
  we then obtain that
  $$\sum\Sb I\in\Cal I\\p\mid P([\![1\in I]\!],\ldots,[\![k\in I]\!])\endSb
  e^{2\pi i\sum_{s\in I}\mu_s}\eq0\ (\mo\ p),$$
  i.e., $|\{I\in\Cal I:\, \sum_{s\in I}\mu_s=\sum_{s\in J}\mu_s\}|
  e^{2\pi i\sum_{s\in J}\mu_s}\eq0\ (\mo\ p)$,
  which is impossible since $p>|\Cal I|$. This concludes our proof. \qed

\Remark\ 3.3. Clearly Corollary 3.2 implies [S99, Theorem 1(i)].
\medskip

 \heading{4. Proof of Theorem 3.1 and a Characterization of $m$-Covers of $\Z$}\endheading

 At first we introduce some notations.
 For a real number $\al$, we let $\{\al\}$ denote the fractional part of $\al$.
 For a polynomial $f(x_1,\ldots,x_k)$ over the field $\C$ of complex numbers,
 we use $[x_1^{j_1}\cdots x_k^{j_k}]f(x_1,\ldots,x_k)$ to represent
 the coefficient of the monomial $x_1^{j_1}\cdots x_k^{j_k}$ in $f(x_1,\ldots,x_k)$.
Also, we fix a finite system (1.2) of residue classes, and set
 $I_z=\{1\ls s\ls k:\, z\in a_s(n_s)\}$ for $z\in\Z$. Note that $|I_z|=w_A(z)\gs m(A)$ for all $z\in\Z$.

 \proclaim{Lemma 4.1} Let $A=\{a_s(n_s)\}_{s=1}^k$
 and let $f(x_1,\ldots,x_k)\in\C[x_1,\ldots,x_k]$ with $\deg f\ls
m(A)$. Let $m_1,\ldots,m_k\in\Z$. If $[\prod_{s\in I_z}x_s]f(x_1,\ldots,x_k)=0$ for all
$z\in\Z$, then we have $\psi(\theta)=0$ for any $0\ls \theta<1$,
where
$$\psi(\theta)=\sum\Sb I\se[1,k]\\\{\sum_{s\in I}m_s/n_s\}=\theta\endSb
(-1)^{|I|}f([\![1\in I]\!],\ldots,[\![k\in I]\!])e^{2\pi i\sum_{s\in I}a_sm_s/n_s}.$$
The converse holds when
$m_1,\ldots,m_k$ are relatively prime to $n_1,\ldots,n_k$
respectively.
\endproclaim

\Proof. Let
$$S=\bg\{\bg\{\sum_{s\in I}\f{m_s}{n_s}\bg\}:\, I\se[1,k]\bg\}.\tag4.1$$
By [S07, Lemma 1], for any $z\in\Z$ we have
$$\sum_{\theta\in S}e^{-2\pi iz\theta}\psi(\theta)
=(-1)^kc(I_z)\prod^k\Sb s=1\\s\not\in I_z\endSb\l(e^{2\pi i(a_s-z)m_s/n_s}-1\r),\tag4.2$$
where $c(I_z)=[\prod_{s\in I_z}x_s]f(x_1,\ldots,x_k)$.

Observe that we must have $c(I_z)=0$
if $\psi(\theta)=0$ for all $0\ls\theta<1$
and each $m_s$ is relatively prime to $n_s$.

Assume that $c(I_z)=0$ for all $z\in[a,\,a+|S|-1]$ where $a\in\Z$.
Then $\sum_{\theta\in S}e^{-2\pi i n\theta}(e^{-2\pi ia\theta}\psi(\theta))=0$ for every
$n\in[0,|S|-1]$. Note that the Vandermonde determinant
$|(e^{-2\pi i\theta})^n|_{n\in[0,|S|-1],\ \theta\in S}$ does not
vanish. So $\psi(\theta)=0$ for each $\theta\in S$. If
$0\ls\theta<1$ and $\theta\not\in S$, then we obviously have
$\psi(\theta)=0$.

 The proof is now complete. \qed

 \Remark\ 4.1. Let $A=\{a_s(n_s)\}_{s=1}^k$ and let $S$ be the set given by (4.1) where
 $m_1,\ldots,m_k\in\Z$ are relatively prime to $n_1,\ldots,n_k$ respectively.
 Let $y$ be any integer with $|I_y|=w_A(y)=m(A)$.
 By Lemma 4.1 and its proof, for any $a\in\Z$ there is a $z\in[a,a+|S|-1]$
 such that $[\prod_{s\in I_z}x_s]\prod_{s\in I_y}x_s\not=0$,
 hence $I_y=I_z$ and $z\in y([n_s]_{s\in I_y})$, where
 $[n_s]_{s\in I_y}$ is the least common multiple of those $n_s$ with $s\in I_y$.
 Therefore $|S|\gs[n_s]_{s\in I_y}$, and we also get
 the following local-global result of Sun [S95, S96, S04]:
  {\it $\{a_s(n_s)\}_{s=1}^k$ forms an $m$-cover of $\Z$ if it covers $|S|$
  consecutive integers at least $m$ times}.
 In the case $m=1$, this local-global principle was conjectured by Erd\H os
 with $|S|$ replaced by $2^k$ (cf. [CV]).
 The reference [S05b] contains a local-global theorem of another type.

\medskip

 Now we use Lemma 4.1 to characterize $m$-covers of $\Z$.
\proclaim{Theorem 4.1} Let $m\in\Z^+$, and let $P_0(x),\ldots,P_{m-1}(x)\in\C[x]$
have degrees $0,\ldots,m-1$ respectively. Let $m_1,\ldots,m_k\in\Z$ and $\mu_1,\ldots,\mu_k\in\C$.
If $A=\{a_s(n_s)\}_{s=1}^k$ forms an $m$-cover of $\Z$, then we have
$$\sum\Sb I\se[1,k]\\\{\sum_{s\in I}m_s/n_s\}=\theta\endSb(-1)^{|I|}P_n\(\sum_{s\in I}\mu_s\)
e^{2\pi i\sum_{s\in I}a_sm_s/n_s}=0\tag4.3$$
for all $0\ls\theta<1$ and $n\in[0,m-1]$. The converse holds provided that
$m_1,\ldots,m_k$ are relatively prime to $n_1,\ldots,n_k$ respectively,
and $\mu_1,\ldots,\mu_k$ are all nonzero.
\endproclaim
\Proof. For $n\in[0,m-1]$ set $f_n(x_1,\ldots,x_k)=P_n(\sum_{s=1}^k\mu_sx_s)$.
If $A=\{a_s(n_s)\}_{s=1}^k$ is an $m$-cover of $\Z$, then $\deg f_n<m\ls m(A)$ and
hence $[\prod_{s\in I_z}x_s]f_n(x_1,\ldots,x_k)=0$ for all $z\in\Z$ since $|I_z|=w_A(z)\gs m$,
therefore (4.3) holds for any $0\ls\theta<1$ in view of Lemma 4.1.

Now assume that $m_1,\ldots,m_k$ are relatively prime to $n_1,\ldots,n_k$
respectively and that $\mu_1\cdots\mu_k\not=0$. Suppose that $m(A)=n\in[0,m-1]$. Then
there is a $z\in\Z$ such that $|I_z|=w_A(z)=n$. As (4.3) holds for all $0\ls\theta<1$, by
Lemma 4.1 the coefficient $c:=[\prod_{s\in I_z}x_s]f_n(x_1,\ldots,x_k)$ vanishes. On the other hand,
$$\align c=&\[\prod_{s\in I_z}x_s\]P_n\(\sum_{s=1}^k\mu_sx_s\)
=[x^n]P_n(x)\times\[\prod_{s\in I_z}x_s\]\(\sum_{s=1}^k\mu_sx_s\)^n
\\=&[x^n]P_n(x)\times\f{n!}{\prod_{s=1}^k[\![s\in I_z]\!]!}\prod_{s\in I_z}\mu_s\not=0.
\endalign$$
This contradiction concludes our proof. \qed

\Remark\ 4.2. In the case $\mu_s=m_s/n_s\ (1\ls s\ls k)$ and $P_n(x)=\bi xn\ (n\in[0,m-1])$,
Theorem 4.1 is equivalent to a characterization
of $m$-covers of $\Z$ obtained by the author [S95, S96] via an analytic method.
\medskip

 Let $a$ be an integer and let $p$ be a prime.
 Fermat's little theorem tells that we can characterize whether
 $p$ divides $a$ as follows:
 $$[\!]p\mid a]\!]\eq1-a^{p-1}\ (\mo\ p).$$
 To handle general abelian $p$-groups in a similar way, we need to
 characterize whether a given power of $p$ divides $a$.
 Thus, our following lemma is of technical importance.
 (It appeared even in the first version of this paper
 posted as {\tt arXiv:math/NT/0305369} on May 24, 2003.)

 \proclaim{Lemma 4.2} Let $p$ be a prime, and let $h\in\N$ and $a\in\Z$.
 Then we have the following congruence
$$\bi{a-1}{p^h-1}\eq[\![p^h\mid a]\!]\ \ (\mo\ p).\tag4.4$$
\endproclaim
\Proof. (4.4) is trivial if $h=0$. Below we let $h>0$.

Write $a=p^hq+r$ where $q,r\in\N$ and $r<p^h$. For $j\in[1,p^h-1]$, if we write
$j=p^{h_j}q_j$ with $h_j\in\N$, $q_j\in\Z^+$ and $p\nmid q_j$, then $0\ls h_j<h$
and
$$\f{p^hq\pm j}{p^h\pm j}=\f{p^{h-h_j}q\pm q_j}{p^{h-h_j}\pm q_j}\eq1\ (\mo\ p).$$
Thus, when $r=0$ we have
$$\bi{a-1}{p^h-1}=\prod_{j=1}^{p^h-1}\f{p^hq-j}{p^h-j}\eq1\ (\mo\ p).$$
In the case $0<r\ls p^h-1$,
$$\bi{a-1}{p^h-1}=\prod_{j=1}^{r-1}\f{p^hq+j}j
\times\f{p^hq}r\times\prod_{j=1}^{p^h-1-r}\f{p^hq-j}{p^h-j}
\eq\f{p^hq}r\eq0\ (\mo\ p).$$
Therefore (4.4) holds. \qed

\medskip
\noindent{\tt Proof of Theorem 3.1}. Write
$$\{P([\![1\in I]\!],\ldots,[\![k\in I]\!])\ \mo\ p:\, I\in\Cal I\}=\{r\ \mo\ p:\, r\in R\}$$
with $R\se[0,p-1]$.
If $0\not\in R$, then (3.3) holds
trivially. So we assume $0\in R$ from now on.

Suppose that $G\cong \Z_{p^{h_1}}\oplus\cdots\oplus\Z_{p^{h_l}}$
where $h_1,\ldots,h_l\in\N$. (When $|G|=1$ we have $G\cong\Z_{p^0}$.)
We can identify $c\in G$ with a vector
$$(c^{(1)}\ \mo\ p^{h_1},\ldots,c^{(l)}\ \mo\ p^{h_l}),$$
and identify $c_s\ (s\in [1,k])$ with a vector
$$(c_{s}^{(1)}\ \mo\ p^{h_1},\ldots,c_{s}^{(l)}\ \mo\ p^{h_l}),$$
where $c^{(t)}$ and $c_{s}^{(t)}$ are integers for $t=1,\ldots,l$.

Let $\theta=\{\al\}$ and
$$\align
f(x_1,\ldots,x_k)=&\prod^l_{t=1}\bi{\sum_{s=1}^kc_{s}^{(t)}x_s-c^{(t)}-1}{p^{h_t}-1}
\\&\times\bi{\sum_{s=1}^km_sx_s/n_s-\al-1}{p^h-1}\prod_{r\in R\sm\{0\}}(P(x_1,\ldots,x_k)-r).
\endalign$$
For any $I\se[1,k]$ with $\{\sum_{s\in I}m_s/n_s\}=\theta$, by Lemma 4.2 we have
$$\align f([\![1\in I]\!],\ldots,[\![k\in I]\!])
\eq&[\![I\in\Cal I]\!]\prod_{r\in R\sm\{0\}}(P([\![1\in I]\!],\ldots,[\![k\in I]\!])-r)
\\\eq&[\![I\in \Cal I\ \&\ p\mid P([\![1\in I]\!],\ldots,[\![k\in I]\!])]\!]C\ (\mo\ p),
\endalign$$
where $C=\prod_{r\in R\sm\{0\}}(-r)\not\eq0\ (\mo\ p)$.
Thus
$$\sigma:=\sum\Sb I\se[1,k]\\\{\sum_{s\in I}m_s/n_s\}=\theta\endSb
(-1)^{|I|}f([\![1\in I]\!],\ldots,[\![k\in I]\!])e^{2\pi i\sum_{s\in I}a_sm_s/n_s}$$
is congruent to the left-hand side of (3.3) times $C\in\Z\sm p\Z$ modulo $p$.

Suppose that
$(3.2)$ fails. Then
$$\deg f\ls \sum_{t=1}^l(p^{h_t}-1)+p^h-1+(|R|-1)d\ls d^*(G)+p^h-1+\Delta<m(A),$$
and hence $\sigma=0$
in light of Lemma 4.1.
Therefore (3.3) holds. We are done. \qed

\enddocument